\documentclass[11pt]{article}
\usepackage{amssymb,amsfonts,amsmath, psfrag,eepic,colordvi}
\newtheorem{theo}{Theorem}

\makeatletter \@addtoreset{equation}{section}
\@addtoreset{theo}{}\makeatother

\def\qed{\hfill \rule{4pt}{7pt}}
\def\pf{\noindent {\it Proof.} }

\def\pf{\noindent {\it Proof.} }

\def\qed{\hfill \rule{4pt}{7pt}}

\begin{document}
\begin{center}
{\LARGE\bf
 Bijective Proofs of  Identities  from Colored Binary Trees}\\[16pt]

 Sherry H.F. Yan\footnote{Research supported by NSFC.}
 \\[6pt]
Department of Mathematics, Zhejiang Normal University\\
\small Jinhua 321004, P.R. China\\[5pt]
hfy@zjnu.cn
\end{center}

\noindent {\bf Abstract. } This note  provide bijective proofs of
two combinatorial identities involving  generalized Catalan number
$C_{m,5}(n)={m\over 5n+m}{5n+m\choose n}$ recently proposed by Sun.

\vskip 8mm



\vskip 1cm

\section{Introduction  }
Recently,  by using generating functions and Lagrange inversion
formula, Sun \cite{Sun}   deduced the following identity involving
generalized Catalan number $C_{m,5}(n)={m\over 5n+m}{5n+m\choose
n}$, i.e.,
\begin{equation}\label{eq.1}
\begin{array}{ll}
 &\sum_{p=0}^{[n/ 4]} {m\over 5p+m}{5p+m\choose
p}{n+p+m-1\choose n-4p} \\
&\\
 &= \sum_{p=0}^{[n/2]}(-1)^{p}{m\over m+n}{m+n+p-1\choose
p}{m+2n-2p-1\choose n-2p},
\end{array}
\end{equation}
which, in the case $m=1$, leads to
\begin{equation}\label{eq.2}
  \sum_{p=0}^{[n/ 4]} {1\over 5p+1}{5p+m\choose
p}{n+p \choose n-4p}= \sum_{p=0}^{[n/2]}(-1)^{p}{1\over
n+1}{n+p\choose p}{2n-2p\choose n-2p}. \end{equation}
 \vskip 5mm
In this note, we give a parity reversing  involution on colored
binary trees which leads to
 a combinatorial interpretation of Formula (\ref{eq.2}).  We make a
 simple variation of the bijection between colored ternary trees and
 binary trees proposed by Sun \cite{Sun} and
 find a correspondence between certain class of  binary trees and the set of colored 5-ary trees. The generalization
 of
 the parity reversing involution and the bijection to forests of
 colored binary trees and forests of colored 5-ary trees leads to a
 bijective proof of Formula (\ref{eq.1}).

 \section{ A parity reversing involution on colored binary trees}
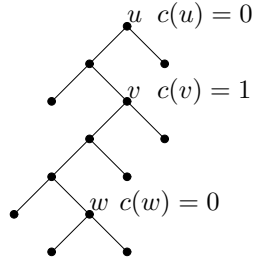
\begin{figure}[h,t]
\begin{center}
\begin{picture}(50,120)
\setlength{\unitlength}{5mm} \linethickness{0.4pt}

\put(2,7){\circle*{0.2}}
 \put(2,7){\line(-1,-1){1}}\put(1,6){\circle*{0.2}}
\put(0,5){\circle*{0.2}} \put(1,6){\line(-1,-1){1}}
 \put(1,6){\line(1,-1){1}}\put(2,5){\circle*{0.2}}
 \put(2,7){\line(1,-1){1}}\put(3,6){\circle*{0.2}}
\put(2,5){\line(-1,-1){1}}\put(1,4){\circle*{0.2}}
\put(2,5){\line(1,-1){1}}\put(3,4){\circle*{0.2}}
\put(1,4){\line(-1,-1){1}}\put(0,3){\circle*{0.2}}
\put(1,4){\line(1,-1){1}}\put(2,3){\circle*{0.2}}
\put(0,3){\line(-1,-1){1}}\put(-1,2){\circle*{0.2}}
\put(0,3){\line(1,-1){1}}\put(1,2){\circle*{0.2}} \put(2,7.1){\small
$u$}\put(2.8,7.1){\small $c(u)=0$} \put(2,5.1){\small
$v$}\put(2.8,5.1){\small $c(v)=1$}
\put(1,2){\line(-1,-1){1}}\put(0,1){\circle*{0.2}}
\put(1,2){\line(1,-1){1}}\put(2,1){\circle*{0.2}} \put(1,2.1){\small
$w$}\put(1.8,2.1){\small $c(w)=0$}
\end{picture}
\end{center}
\caption{ A colored binary tree $B$ with $c(B)=1$.} \label{fig2}
\end{figure}

 In this section, we provide a parity reversing involution on a
 class of colored binary trees. Before introducing the involution,
 we recall
 some definitions and notations. Let $\mathcal{B}_n$
 denote the set of complete binary trees with $n$ internal vertices.
 Let $B\in \mathcal{B}_n$ and $P=v_0v_1\ldots v_k$ be a path of
 length $k$ of $B$ (viewing from the root). $P$ is called a {\em
 L-path} if $v_i$ is a left child of $v_{i-1}$ for $1\leq i\leq k$.
 $P$ is called a {\em maximal} L-path if there exists no vertex such
 that $uP$ or $Pu$ forms a L-path. Suppose that $P=v_0v_1 \ldots
 v_k$ is a maximal L-Path, then  $v_0$ is
 called an {\em   initial} vertex of $B$ and $P$ is called the {\em associate }  path of
 $v_0$. Denote by $l(v)$ the length of the associate  maximal L-path of $v$.
  A {\em colored binary tree} is a binary tree in which each initial vertex $v$ is assigned a color $c(v)$
such that $0\leq c(v)\leq [l(v)/2]$.
 The {\em color number} of a
colored binary tree $B$, denoted by $c(B)$, is equal to the sum of
all the colors of initial vertices of $B$. A colored binary tree $B$
with $c(B)=1$ is illustrated in Figure \ref{fig2}.

 Let
$\mathcal{CB}_{n,p}=\{B| B\in B_n \, {\mbox with}\,  c(B)=p\}$ and
$f_{n,p}$ be its cardinality.
    Define  $\mathcal{CB}_{n
 }=\bigcup_{p=0}^{[n/2]}\mathcal{CB}_{n,p}$.
 Let $f(x,y)=\sum_{n\geq
 0}\sum_{p=0}^{[n/2]}f_{n,p}y^px^n$ be the ordinary generating function for $f_{n,p}$ with the assumption $f(0,0)=1$.
Given a colored binary tree $B\in \mathcal{CB}_{n}$, let $v$ be the
root of $B$. Suppose that $c(v)=k$, then  $l(v)\geq 2k$.  Then the
generating function for the number of colored binary trees whose
root $v$ colored by $k$  is equal to $\sum_{k\geq
0}y^k(xf(x,y))^{2k}(xf(x,y))^{l(v)-2k}$.   Summing over all the
possibilities for $k\geq0$, we arrive at
 $$
 f(x,y)={1 \over {(1-yx^2f^2(x,y))}(1-xf(x,y))}.
 $$
 By applying Lagrange inversion formula, we get
 $$
 f^m(x,y)=\sum_{n\geq 0}\sum_{p=0}^{[n/2]}{m\over
n+m} {m+n+p-1\choose p}{m+2n-2p-1\choose n-2p}y^px^n,
 $$
which in the case of $m=1$, reduces to
$$
f(x,y)=\sum_{n\geq 0}\sum_{p=0}^{[n/2]}{1\over n+1}{n+p\choose p} {
2n-2p \choose n-2p}y^px^n.
$$

  Let $B\in \mathcal{CB}_{n
 }$ and $v$ be an initial vertex of $B$. If $c(v)=2k$ and $l(v)=4k $ or $4k+1$
 for some nonnegative integer $k$, then we say $v$ is a {\em proper}  initial vertex,
 otherwise it is said to be {\em improper}. A colored  binary
 tree $B$ is called {\em proper} if all the initial vertices are proper; otherwise, $B$ is said to be {\em improper}.
 Denote by $\mathcal{CB}'_{n}$
the set of all proper   binary trees with $n$ internal vertices.
 Define $\mathcal{ECB}_n$ and $\mathcal{OCB}_n$ to be the sets of colored
 binary trees with $n$ internal vertices whose color numbers are even
 and odd, respectively.
\begin{theo}\label{theo1}
There is a parity reversing involution $\phi$ on the set of improper
colored binary trees with $n$ internal vertices. Furthermore,
$$
|\mathcal{ECB}_n|-|\mathcal{OCB}_n|=|\mathcal{CB}'_{n}|.
$$
\end{theo}
\pf Let $B$ be an improper colored binary trees with $n$ internal
vertices with $c(B)=p$. Traverse the binary tree $B$ by
 depth first search, let $v$ be the first improper vertex
 traversed. If $c(v)$ is even,  then let $\phi(B)$ be a colored
 binary tree obtained from $B$ by  coloring $v$ by $c(v)+1$.
   If  $c(v)$ is odd, then let $\phi(B)$ be a colored binary tree
obtained from  $B$ by coloring $v$ by $c(v)-1$. Obviously,  the
obtained binary tree $\phi(B)$ are improper colored binary trees in
both cases.  Furthermore, in the former case the color number of
$\phi(B)$ is $p+1$, while in the latter case, the color number is
$p-1$. Hence $\phi$ is an involution on the set of improper colored
binary trees with $n$ internal vertices, which reverses the parity
of the color number. Since each initial vertex $v$ of a proper
binary tree is colored by an even number, it is clear that
$\mathcal{CB}'_{n}\subseteq \mathcal{ECB}_n$. Hence, we have
$$
|\mathcal{ECB}_n|-|\mathcal{OCB}_n|=\mathcal{CB}'_{n}.
$$ \qed

From Theorem \ref{theo1}, we see that the right side summand of
Formula (\ref{eq.2}) counts the number of proper colored binary
trees, that is,
$$
|\mathcal{CB}'_n|=\sum_{p=0}^{[n/2]}(-1)^{p}{1\over n+1}{n+p\choose
p}{2n-2p\choose n-2p}.
$$

\section{The bijective proof }
 A (complete) {\em $k$-ary} tree is an ordered tree in which each
 internal vertex has $k$ children. The number of $k$-ary trees with
 $n$ internal vertices is counted by generalized Catalan number $C(1,k)(n)={1\over kn+1}{kn+1\choose
 n}$ \cite{st}.  A {\em colored} $5$-ary tree is a $5$-ary tree
  in which each vertex is assigned a nonnegative integer called {\em color
 number}, denoted by $c_v$. Let $\mathcal{T}_{n,p}$ denote the set of colored $5$-ary
 trees with $p$ internal vertices  such that the sum of all the color numbers of each tree is
 $n-4p$. Denote by $\mathcal{T}_n=\bigcup_{p=0}^{[n/ 4]} \mathcal{T}_{n,p}$.  Let $B\in \mathcal{CB}'_n$ be a proper binary tree with $n$
 internal vertices.  Since each initial vertex $v$ of $B$ is colored by $[l(v)/2]$,  we can discard all the colors of vertices of $B$.

\begin{figure}[h,t]
\begin{center}
\begin{picture}(300,100)
\setlength{\unitlength}{5mm} \linethickness{0.4pt}

\put(0,5){\circle*{0.2}}
 \put(0,5){\line(2,1){2}}\put(0, 4.5){\small$T_1$}
\put(2,6){\circle*{0.2}} \put(2,6){\line(0,1){1}}
\put(2,7){\circle*{0.2}} \put(1.6, 6){\small$v$} \put(2.4,
6){\small$c_v=2$} \put(2,6){\line(-1,-1){1}}\put(1,5){\circle*{0.2}}
\put(1, 4.5){\small$T_2$}
\put(2,6){\line(1,-1){1}}\put(3,5){\circle*{0.2}} \put(3,
4.5){\small$T_4$} \put(2,6){\line(2,-1){2}}\put(4,5){\circle*{0.2}}
\put(4,
4.5){\small$T_5$}\put(2,6){\line(0,-1){1}}\put(2,5){\circle*{0.2}}
\put(2, 4.5){\small$T_3$}
\put(6,6){$\Leftrightarrow$} \put(6 , 0){$(a)$}
\put(8,7){\circle*{0.2}}\put(8,7){\line(0,-1){1}}\put(8,6){\circle{0.2}}

\put(8,6){\line(-1,-1){1}}\put(7,5){\circle{0.2}}
\put(8,6){\line(1,-1){1}}\put(9,5){\circle{0.2}}
\put(9,5){\line(-1,-1){1}}\put(8,4){\circle{0.2}}
\put(9,5){\line(1,-1){1}}\put(10,4){\circle{0.2}}\put(9.5,
3.9){\small$v$}

 \put(10,4){\line(-2,-1){2}}\put(8,3){\circle{0.2}}
\put(8, 2.5){\small$T_1$}
\put(10,4){\line(2,-1){2}}\put(12,3){\circle{0.2}} \put(12,
2.5){\small$T_5$} \put(10,4){\line(-1,-1){1}}\put(9,3){\circle{0.2}}
\put(9, 2.5){\small$T_2$}
\put(10,4){\line(1,-1){1}}\put(11,3){\circle{0.2}} \put(11,
2.5){\small$T_4$} \put(10,4){\line(0,-1){1}}\put(10,3){\circle{0.2}}
\put(10, 2.5){\small$T_3$}
\put(14,6){$\Leftrightarrow$} \put(13.5 , 0){$(b)$}
\put(18,7){\circle*{0.2}}\put(18,7){\line(0,-1){1}}\put(18,6){\circle{0.2}}

\put(18,6){\line(-1,-1){1}}\put(17,5){\circle{0.2}}
\put(18,6){\line(1,-1){1}}\put(19,5){\circle{0.2}}
\put(19,5){\line(-1,-1){1}}\put(18,4){\circle{0.2}}
\put(19,5){\line(1,-1){1}}\put(20,4){\circle{0.2}}\put(19.5,
3.9){\small$v$}

 \put(20,4){\line(-1,-1){1}}\put(19,3){\circle{0.2}}\put(18.3,
2.9){\small$v_1$}
\put(20,4){\line(1,-1){1}}\put(21,3){\circle{0.2}}\put(21.3,
2.9){\small$T_5$}

\put(19,3){\line(-1,-1){1}}\put(18,2){\circle{0.2}}\put(17.3,
1.9){\small$v_2$}
\put(19,3){\line(1,-1){1}}\put(20,2){\circle{0.2}}\put(20.3,
1.9){\small$T_4$}

\put(18,2){\line(-1,-1){1}}\put(17,1){\circle{0.2}}\put(16.3,
0.9){\small$v_3$}
\put(18,2){\line(1,-1){1}}\put(19,1){\circle{0.2}}\put(19.3,
0.9){\small$T_3$}

\put(17,1){\line(-1,-1){1}}\put(16,0){\circle{0.2}}\put(15.2,
-0.1){\small$T_1$}
\put(17,1){\line(1,-1){1}}\put(18,0){\circle{0.2}}\put(18.3,
-0.1){\small$T_2$}
\end{picture}
\end{center}
\caption{The bijection $\sigma$.} \label{fig1}
\end{figure}
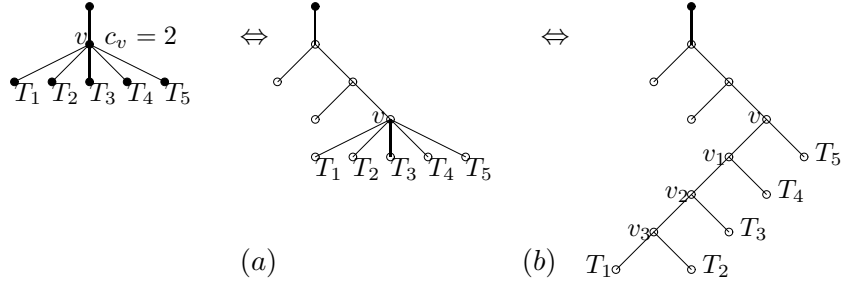

  Now we construct a map
 $\sigma$ from $\mathcal{T}_n$ and $
 \mathcal{CB}'_n$ as follows:
 \begin{itemize}
 \item[Step 1.] For each vertex $v$ of $T\in \mathcal{T}_n$ with
 color number $c_v=k$, remove the color number and add a  path $P=v_1v_2\ldots
 v_k$   to $v$ such that $v$ is a right child of
 $v_k$ and $v_1$ is a child of the father of $v$, and annex a left
 leaf to $v_i$ for $1\leq i\leq k$. See Figure \ref{fig1}(a) for
 example.

 \item [Step 2.]  Let $T^*$ be the tree obtained from $T$ by Step
 $1$. For any internal vertex $v$ of $T^*$ which has $5$ children, let $T_1, T_2, T_3, T_4,
 T_5$ be the five subtrees of $v$.  Annex a  path $P=v_1v_2v_3$ to $v$ such that $v_1$ is a left child of
 $v$, then
 take $T_1$ and $T_2$ as the left and right subtree of $v_3$, take
 $T_3$ as the right subtree of $v_2$, and take $T_4$ as the subtree
 of $v_1$.  See Figure \ref{fig1}(b) for
 example.
 \end{itemize}

It is clear that the obtained tree $\sigma(T)$ is a proper binary
tree with $n$ internal vertices.

Conversely, we can obtain a colored 5-ary tree from a proper binary
tree by a similar   procedure from binary trees to colored ternary
trees given by Sun \cite{Sun}. We omit the reverse map of $\sigma$
here.

 \begin{theo}\label{theo2}
  The map $\sigma$ is a bijection between $\mathcal{T}_n$ and $
 \mathcal{CB}'_n$.
 \end{theo}

 Given a   5-ary tree $T$ with $p$ internal vertices,there are
 totally $5p+1$ vertices. Choose $n-4p$ vertices repeatedly, and
 define the color of the vertex as the number of times being
 chosen. Then there are ${n+p\choose n-4p}$ colored 5-ary trees in
 $\mathcal{T}_n$ obtained from $T$. Since there are ${1\over 5p+1}{5p+1\choose p}$
 5-ary trees with $p$ internal vertices, the involution $\phi$ and
 the bijection $\sigma$ lead to a combinatorial proof of Formula
 (\ref{eq.2}).

In order to prove Formula (\ref{eq.1}), we consider   the forest of
colored binary trees $F=(B_1, B_2, \ldots, B_m)$ with $n$ internal
vertices and $m$ components where $B_i\in \mathcal{CB}_{n_i}$ and
$n_1+n_2+\ldots+n_m=n$.  Define the {\em color number} of $F$ as the
sum of all the color numbers of $B_i$, where $1\leq i\leq m$. It is
easy to check that the number of forests of colored binary trees
with $m$ components and $n$ internal vertices, whose color number
equals $p$, is counted by

 $$[y^px^n]f^m(x,y)= {m\over
n+m} {m+n+p-1\choose p}{m+2n-2p-1\choose n-2p}.$$

$F$ is said to be a {\em proper} forest if each $B_i$ is proper;
otherwise, $F$ is said to be {\em improper}. Now we can modify the
involution $\phi$ as follows:  suppose that $B_k$ be the leftmost
improper   binary tree, then let $\phi(F)$ be the forest of colored
binary tree obtained from $F$ by changing $B_k$ to $\phi(B_k)$. From
Theorem \ref{theo1}, we see $\phi(F)$ is an improper forest of
colored binary trees with $n$ internal vertices and $m$ components.
Hence the modified involution $\phi$ is an involution on the set of
improper forests of  colored  binary trees, which reverses the
parity of the color numbers of the forests. Hence the right side
summand of Formula (\ref{eq.1}) counts the number of proper forests
of colored binary trees with $n$ internal vertices and $m$
components.

Let $F=(T_1, T_2, \ldots, T_m)$ be a forest of 5-ary trees such that
$T_i\in \mathcal{T}_{n_i}$ and $n_1+n_2+\ldots+n_m=p$.
 Define $\sigma(F)=(\sigma(T_1), \sigma(T_2),
\ldots, \sigma(T_m))$.  From Theorem \ref{theo2}, it is clear that
$\sigma(F)$ is a proper forest  of   binary trees.  Note that there
are totally $m+5p$ vertices in a forest  $F $ of 5-ary trees with
$p$ internal vertices and $m$ components, so there are
${n+p+m-1\choose n-4p}$ forests of colored 5-ary trees with $m$
components, $p$ internal vertices and the sum of the color numbers
equal to $n-4p$.  Hence the modified involution $\phi$ and the
modified bijection $\sigma$ leads to a bijective proof of Formula
(\ref{eq.1}).

\small

\end{document}